\documentclass[twoside]{amsart}

\usepackage{latexsym, amsmath, amssymb, amsfonts, amscd, amsthm}

\newcommand{\Ce}{{\mathbb C}}
\renewcommand{\Re}{{\mathbb R}}

\newcommand{\Ne}{{\mathbb N}}

\newcommand{\pa}{\parallel}

\newcommand{\ie}{{\it i.e. }}

\theoremstyle{plain}
\newtheorem{theorem}{Theorem}[section]
\newtheorem{lemma}[theorem]{Lemma}
\newtheorem{proposition}[theorem]{Proposition}
\newtheorem{corollary}[theorem]{Corollary}

\theoremstyle{definition}
\newtheorem{remark}[theorem]{Remark}
\newtheorem{example}[theorem]{Example}

\newtheorem{definition}[theorem]{Definition}
\newtheorem{conventions}[theorem]{Conventions}
\newtheorem*{acknowledgments}{Acknowledgments}

\begin{document}

\title{Strict Quantizations of Almost Poisson Manifolds}
\author{Hanfeng Li}

\address{Department of Mathematics \\
University of Toronto \\
Toronto ON M5S 3G3, CANADA}
\email{hli@fields.toronto.edu}
\date{May 18, 2003}

\begin{abstract}
We show the existence of (non-Hermitian) strict quantization for
every almost Poisson manifold.
\end{abstract}

\maketitle

\section{Introduction}
\label{intro:sec}

In the passage from classical mechanics to quantum mechanics, smooth functions
on symplectic manifolds (more generally, Poisson manifolds) are replaced by operators
on Hilbert spaces, and the Poisson bracket of smooth functions are replaced by commutators of operators.
When one thinks of classical mechanics as limits of quantum mechanics, the Poisson bracket becomes
limits of commutators. Based on the general theory of formal deformations of algebras \cite{Gerstenhaber64},
F. Bayer et al. \cite{BFFLS78} initiated the study of {\it deformation quantization} of Poisson manifolds.

Let $M$ be a Poisson manifold. Denote $C^{\infty}(M)$ the space of
smooth $\Ce$-valued functions on $M$, and denote
$C^{\infty}(M)[[\hbar]]$ the space of formal power series with
coefficients in $C^{\infty}(M)$. Recall that a {\it star product}
on $M$ is a $\Ce[[\hbar]]$-bilinear associative multiplication $*$
on $C^{\infty}(M)[[\hbar]]$ of the form
\begin{eqnarray*}
f*g=\sum^{\infty}_{r=0}C_r(f,g)\hbar^r,\, \, \mbox{ for } f, g\in C^{\infty}(M),
\end{eqnarray*}
where $C_0(f,g)=fg, \, f*g-g*f\equiv \{f,g\}i\hbar \mod \hbar^2$,
and each $C_r(\cdot, \cdot)$ is a bidifferential operator. The
algebra $(C^{\infty}(M)[[\hbar]], *)$ is called a {\it deformation
quantization} of $M$. The existence of deformation quantizations
for any symplectic manifold was proven first by De Wilde and
Lecomte \cite{DL83}. The general case of Poisson manifolds was
proven by Kontsevich \cite{Kontsevich97}.

In deformation quantizations $\hbar$ is only a formal parameter,
and elements in $C^{\infty}(M)[[\hbar]]$ are not operators on
Hilbert spaces. In order to study quantizations in a stricter
sense, Rieffel introduced {\it strict deformation quantization} of
Poisson manifolds \cite{Rieffel89, Rieffel91}, and showed that
noncommutative tori arise naturally as strict deformation
quantizations of certain Poisson brackets on the ordinary torus.
Later, Landsman introduced a weaker notion {\it strict
quantization} to accommodate some other interesting examples such
as Berezin-Toeplitz quantization of K$\ddot{a}$hler manifolds.
Recall the definition of strict quantization as formulated in
\cite{Rieffel98a, Landsman98b}:
\begin{definition} \label{sq:def}
Let $M$ be a Poisson manifold,
and let $C_{\infty}(M)$ be the algebra of continuous functions on $M$ vanishing at $\infty$.
By a \emph{strict
quantization} of M we mean a dense $*$-subalgebra $A$ of
$C_{\infty}(M)$ closed under the Poisson bracket, together with a
continuous field of $C^*$-algebras $\mathcal{A}_{\hbar}$  over a closed
subset $I$ of the real line containing $0$ as a non-isolated
point, and linear maps $\pi_{\hbar}: A\to \mathcal{A}_{\hbar}$ for each
$\hbar\in I$, such that

(1)  $\mathcal{A}_0=C_{\infty}(M)$
 and $\pi_0$ is the canonical inclusion of
     $A$ into $C_{\infty}(M)$,

(2)  the section $(\pi_{\hbar}(f))$ is continuous for every $f\in
A$,

(3)  for all $f,g\in A$ we have
\begin{eqnarray*}
   \lim_{\hbar\to 0}\pa
   [\pi_{\hbar}(f), \, \pi_{\hbar}(g)]/(i\hbar)-\pi_{\hbar}(\{f,g\})\pa
   = 0.
\end{eqnarray*}

If each $\pi_{\hbar}$ is injective, we say that the strict
quantization is \emph{faithful}.  If $A\supseteq C^{\infty}_c(M)$,
the space of compactly supported smooth functions on $M$,
we also say that the strict quantization is \emph{flabby}.
If $(\pi_{\hbar}(f))^*=\pi_{\hbar}(f^*)$ for all $\hbar\in I$ and $f\in A$, we say that
the strict quantization is \emph{Hermitian}.
When a Lie group has an smooth action $\alpha$ on $M$ preserving the
Poisson bracket, if $G$ also has an continuous action $\beta_{\hbar}$
on each $\mathcal{A}_{\hbar}$ such that $\beta_0=\alpha^*$ and
the maps $\pi_{\hbar}$ are all $G$-equivariant, we say that the strict quantization is \emph{$G$-equivariant}.
When
the strict quantization is faithful and $\pi_{\hbar}(A)$ is a
$*$-subalgebra of $\mathcal{A}_{\hbar}$ for each $\hbar$,  it's called a
\emph{strict deformation quantization}.
\end{definition}

Strict quantizations have been constructed for several classes of
Poisson manifolds such as Poisson manifolds coming from actions of
$\Re^d$ \cite{Rieffel93}, quantizable compact K$\ddot{a}$hler
manifolds \cite{Schlichenmaier94}, dual of integrable Lie
algebroids \cite{Rieffel90a, LR01}, compact Riemannian surfaces of
genus $\ge 2$ \cite{KL97, Natsume01, NN01} etc. These
constructions are all global, and the resulting strict
quantizations are Hermitian. However, the progress of the study of
strict quantizations is much slower compared with that of
deformation quantizations--so far there is even no existence
result for general symplectic manifolds. Recently Natsume et al.
\cite{NNP02} constructed strict quantizations for every compact
symplectic manifold $M$ satisfying the topological conditions that
$\pi_1(M)$ is exact and $\pi_2(M)=0$. Roughly speaking, they use
partition of unity to reduce $M$ to Darboux charts, where they can
use the Moyal-Weyl product. Thus their construction is local. It
turns out that the resulting strict quantizations are not
Hermitian.

Recall that an {\it almost Poisson manifold} is a smooth
manifold $M$ equipped with some $\Pi\in \Gamma(\wedge^2TM)$ \cite{CW01}. In
this case, we can still define a bracket $\{f, g\}=\Pi(df, dg)$
for $f, g\in C^{\infty}(M)$,  which is bilinear and
skew-symmetric, and satisfies the Leibniz rule.  And the bracket
satisfies the Jacobi identity if and only if $M$ equipped with
this bracket is actually a Poisson manifold. Clearly we can also
talk about strict quantizations of almost Poisson manifolds.
The main result of this paper is the following:

\begin{theorem} \label{sq for almost Poisson manifold 2:thm}
Let $(M, \Pi)$ be an almost Poisson manifold, and let $\alpha$ be a smooth
action of a Lie group $G$ on $M$ preserving the bracket. If $M$
has a $G$-invariant Riemannian metric, then $M$ has a
$G$-equivariant faithful flabby strict quantization over $I=[0,1]$
with $A=C^{\infty}_c(M)$. In particular, taking $G=\{e\}$ we see
that $M$ has a faithful flabby strict quantization over $[0,1]$
with $A=C^{\infty}_c(M)$.
\end{theorem}

Our construction is also local, but different from the one in
\cite{NNP02}. Actually we shall construct a locally trivial
$C^*$-algebra bundle over $M$ in a canonical way, thus don't need
local charts and partition of unity. But our strict quantizations
are not Hermitian either.

This paper is organized as follows. Though our construction for
strict quantizations of almost Poisson manifolds is only slightly
more complicated for that of symplectic manifolds, the idea is
most natural in the case of symplectic manifolds. So we prove
Theorem~\ref{sq for almost Poisson manifold 2:thm} for symplectic
manifolds first in Section~\ref{sq of SM:sec}. Then we prove
Theorem~\ref{sq for almost Poisson manifold 2:thm} for the general
case in Section~\ref{sq alm P mani:sec}. Our construction depends
on the choice of an inner product on the vector bundle $T^*M\oplus
TM$. We define {\it homotopy} of strict quantizations in
Section~\ref{homotopy:sec}, and show that the homotopy class of
our strict quantizations doesn't depend on the choice of the inner
products. In Section~\ref{local strict quantization:sec} we define
{\it local strict quantizations}, and show that they can't be
Hermitian. We also show that our strict quantizations can't be
restricted to a $*$-subalgebra of $A$ to get a strict deformation
quantization of $M$ unless $\Pi=0$. This gives a negative answer
to a question of Rieffel \cite[Question 25]{Rieffel98a}. In
Section~\ref{functor:sec} we discuss certain functorial properties
of our construction. All of our construction is based on the
existence of {\it asymptotic representations of Heisenberg
commutation relations} (Definition~\ref{Heisenberg:def}). We prove
the existence of such  asymptotic representations in
Section~\ref{asym rep:sec}.

\begin{acknowledgments}
  I would like to thank Marc Rieffel for many helpful discussions and suggestions.
  I also thank Henrique Bursztyn for valuable discussions about deformation quantizations.
\end{acknowledgments}

\section{Strict Quantizations of Symplectic Manifolds}
\label{sq of SM:sec}
Throughout this paper, for a continuous field $D$ of $C^*$-algebras
$\{D_{\hbar,x}\}_{x\in X}$ over a locally compact Hausdorff space $X$ we denote
$\Gamma_{\infty}(D)$ the algebra of continuous sections of $D$ vanishing at $\infty$ \cite{Dixmier77}.

We show the main idea of our construction first. Let $(M, \omega)$
be a symplectic manifold, and let $\Pi\in \Gamma(\wedge^2TM)$ be
the corresponding bivector field as usual. Let $f, g\in
C^{\infty}_c(M)$. Since we think of $\mathcal{A}_{\hbar}$ as
deformations of $C_{\infty}(M)$, we would like to write
$\pi_{\hbar}(f)$ as $f+\tau_{\hbar}(f)$, which makes sense when
$\mathcal{A}_{\hbar}$ contains $C_{\infty}(M)$ as a
$C^*$-subalgebra, and assume that
\begin{eqnarray} \label{vanish:eq}
\pa \tau_{\hbar}(f)\pa \rightarrow 0
\end{eqnarray}
as $\hbar\rightarrow 0$. Assume further that $C_{\infty}(M)$ lies
in the center of $\mathcal{A}_{\hbar}$. Then
\begin{eqnarray*}
[\pi_{\hbar}(f),\pi_{\hbar}(g)]/(i\hbar)-\pi_{\hbar}(\{f,g\})=
[\tau_{\hbar}(f),\tau_{\hbar}(g)]/(i\hbar)-\{f,g\}-\tau_{\hbar}(\{f,g\}).
\end{eqnarray*}
Thus the condition (3) in Definition~\ref{sq:def} becomes $\pa
[\tau_{\hbar}(f),\tau_{\hbar}(g)]/(i\hbar)-\{f,g\}\pa\rightarrow
0$. Notice that $\{f, g\}$ doesn't depend on $f$ and $g$, but
depends only on $df$ and $dg$. So we would like to assume that
$\tau_{\hbar}(f)$ depends only on $df$ linearly. Then we attempt
to  write $\tau_{\hbar}(f)$ as
$\hbar^{\frac{1}{2}}\varphi_{\hbar}(df)$, where
$\varphi_{\hbar}:\Gamma(T^*M)\rightarrow \mathcal{A}_{\hbar}$ is a
linear map. Now $\pa
[\tau_{\hbar}(f),\tau_{\hbar}(g)]/(i\hbar)-\{f,g\}\pa \rightarrow
0$ becomes
\begin{eqnarray} \label{comm conv global:eq}
\pa [\varphi_{\hbar}(df), \varphi_{\hbar}(dg)]-\Pi(df, dg)i\pa
\rightarrow 0.
\end{eqnarray}
We assume further that $\mathcal{A}_{\hbar}=\Gamma_{\infty}(D)$
for some continuous field of $C^*$-algebras $\{D_{\hbar,x}\}_{x\in
M}$ over $M$, and that $D_{\hbar, x}$ contains $T^*M_x$ as a
linear subspace with $\varphi_{\hbar}$ being just pointwise
embedding. Then (\ref{comm conv global:eq}) becomes
\begin{eqnarray} \label{comm conv point:eq}
\pa [u,v]-\Pi_x(u,v)i\pa_{D_{\hbar,x}}\rightarrow 0.
\end{eqnarray}
for all $u,v \in T^*M_x$.

This leads to our definition of asymptotic representation of
Heisenberg commutation relations and that of Heisenberg
$C^*$-algebra $\mathfrak{A}_{2n}$ (which will be $\mathcal{A}_{\hbar, x}$
for $2n=\dim M$) in Definition~\ref{Heisenberg:def} below. In
order to embed $T^*M_x$ into $\mathfrak{A}_{2n}$ without referring
to local basis, we also need an action of the structure group of
$T^*M$ on $\mathfrak{A}_{2n}$. If we consider only $\Pi$, then
the structure group is the symplectic linear group $Sp(2n)$
\cite{Cannas01}, which is too big. By adding a compatible almost
complex structure on $T^*M$ we can reduce the structure group to
the unitary group $U(n)$ (see Lemma~\ref{local unitary
basis:lemma}). Here we recall how $U(n)$ acts on $T^*M_x$.

Let $V$ be a finite dimensional vector space over $\Re$ with a
symplectic structure $\omega$ and a compatible almost complex
structure $J$ \cite{Cannas01}, \ie $J:V\to V$ is linear satisfying
that $J^2=-1$ and $<u,v>:=\omega(u,Jv)$ is an inner product on
$V$. Say $dim V=2n$. Then we can always find basis $u_1,\, \cdots,\,
u_n,\, v_1,\,\cdots,\, v_n$ of $V$ such that under this basis $\omega$
and $J$ have matrix forms $\left(
\begin{array}{cr} 0 & I \\ -I & 0\end{array} \right)$ and $\left(
\begin{array}{cr} 0 & -I \\ I & 0\end{array} \right)$
respectively.  We'll call such basis {\it unitary basis} of $V$.
Notice that a unitary basis is an orthonormal basis under the
induced inner product. If we make $V$ into a complex vector space
by $J$ and identify matrix $X+iY\in M_{n\times n}(\Ce)$, where $X,
Y\in M_{n\times n}(\Re)$, with $\left(
\begin{array}{cr} X & -Y \\ Y & X \end{array} \right) \in
M_{2n\times 2n}(\Re)$, then $U(n)$ is exactly the group of linear
transformations on $V$ taking unitary basis to unitary basis.

\begin{definition} \label{Heisenberg:def}
Let $n\in \Ne$, and let $\Re^{2n}$ be equipped with the standard
symplectic vector space structure $\Omega$ and the standard
compatible almost complex structure, \ie for the standard basis
$e_1, \cdots, e_{2n}$ being a unitary basis.
By an \emph{asymptotic representation of Heisenberg commutation
relations}, we mean a unital $C^*$-algebra $\mathfrak{A}_{2n}$
with a continuous action $\rho$ of $U(n)$ and a $U(n)$-equivariant
$\Re$-linear map $\varphi_{\hbar}:\Re^{2n}\rightarrow
\mathfrak{A}_{2n}$ for each $0<\hbar\le 1$ such that

(1) for any $u,v\in \Re^{2n}$ we have $[\varphi_{\hbar}(u),
\varphi_{\hbar}(v)]\rightarrow \Omega(u,v)i$ as $\hbar\rightarrow
0$;

(2) the map $(0, 1]\rightarrow B(\Re^{2n}, \mathfrak{A}_{2n})$
given by $\hbar\mapsto \varphi_{\hbar}$ is continuous, where \\
$B(\Re^{2n}, \mathfrak{A}_{2n})$ is the Banach space of linear
maps $\Re^{2n}\rightarrow \mathfrak{A}_{2n}$;

(3) $\hbar^{\frac{1}{2}}\pa \varphi_{\hbar}\pa\rightarrow 0$ as
$\hbar\rightarrow 0$;

(4) $\mathfrak{A}_{2n}$ is generated by $\cup_{0<\hbar\le
1}\varphi_{\hbar}(\Re^{2n})$.

The $C^*$-algebra $\mathfrak{A}_{2n}$ will be called a \emph{Heisenberg $C^*$-algebra
of dimension $2n$}.
\end{definition}

\begin{remark} \label{Heisenberg:remark}
The condition (3) is not crucial. Given $\mathfrak{A}_{2n}$ and
$\varphi_{\hbar}$ satisfying the other conditions we can always reparameterize $\varphi_{\hbar}$'s
to make them satisfy (3).
\end{remark}
The main technical part of our construction is the following:

\begin{theorem} \label{asymptotic rep exit:thm}
For each $n\in \Ne$ there exists a Heisenberg $C^*$-algebra
$\mathfrak{A}_{2n}$.
\end{theorem}

Theorem~\ref{asymptotic rep exit:thm} will be proved in
Section~\ref{asym rep:sec}. From now on we'll fix a Heisenberg
$C^*$-algebra $\mathfrak{A}_{2n}$ for each $n\in \Ne$ unless
stated otherwise.

We'll construct a $C^*$-algebra bundle $D$ over $M$ with fibres
all isomorphic to $\mathfrak{A}_{2n}$ and a bundle map
$T^*M\rightarrow D$. Actually we shall do this construction more
generally for almost symplectic bundles, which will be useful in
Section~\ref{sq alm P mani:sec}.

\begin{definition} \label{alm sym bundle:def}
Let $E\rightarrow M$ be a real vector bundle, and let $\Pi \in
\Gamma(\wedge^2E^*)$.
We call the pair $(E, \Pi)$ an \emph{almost Poisson bundle} over
$M$.
If $\Pi$ is nondegenerate everywhere we
call it an \emph{almost symplectic bundle} over
$M$. For two almost Poisson bundles $(E, \Pi_E)$ and $(F, \Pi_F)$ over $M$ we call
a bundle map $\psi:E\rightarrow F$ an \emph{almost Poisson map} if  $\Pi_F(\psi(f),\psi(g))=\Pi_E(f,g)$
for all $f,g\in \Gamma(E)$.
\end{definition}

As the case for symplectic manifolds \cite{Cannas01}, one sees easily that
every almost symplectic bundle $(E, \Pi)$ has compatible almost
complex structures.

\begin{lemma} \label{local unitary basis:lemma}
Let $(E, \Pi)$ be an almost symplectic bundle over $M$ with a
compatible almost complex structure $J$. Say $dim E=2n$. Then for any
$x\in M$ there exists a neighborhood $U$ of $x$  and sections
$f_1, \, \cdots, \, f_n,\, g_1,\, \cdots, \, g_n$ of $E$ on $U$ such that they are a
unitary basis at each point of $U$.
\end{lemma}
\begin{proof}
It is easy to see that we can find sections $f_1, \, \cdots, \, f_n,\, g_1,\, \cdots, \, g_n
$ near $x$ such that under this basis, $\Pi$ is of the
matrix form $\left(
\begin{array}{cr} 0 & I \\ -I &  0\end{array} \right)$ at each
point.  Let $E^{+}=span\{f_1, \cdots, f_n\}$.  Then near $x$ this
is a subbundle of $E$, and $E^{+}_y$ is a Lagrangian subspace of
$E_y$ at each point $y$.  Let $f'_1, \cdots, f'_n$ be $n$ sections
of $E^{+}$ near $x$ such that $f'_1(y), \cdots, f'_n(y)$ is an
orthonormal basis at each point $y$.  Then $f'_1, \,\cdots ,\,  f'_n,\,
Jf'_1,\, \cdots,\,  Jf'_n$ satisfy the requirement.
\end{proof}

Lemma~\ref{local unitary basis:lemma} shows that the structure
group of $(E, \Omega, J)$ is $U(n)$. Let $U(E)$ be the set of
unitary basis of $E$ at all points, then $U(E)$ is a principle
$U(n)$-bundle on $M$. As usual, the action $\rho$ of $U(n)$ on
$\mathfrak{A}_{2n}$ induces a $C^*$-algebra bundle
$D=U(E)\times_{U(n)}\mathcal{A}_{2n}$ over $M$, which is the
quotient of $U(E)\times \mathcal{A}_{2n}$ by the relation $(a,
T)\sim(ag, g^{-1}T)$ with $g\in U(n)$. Then $D$ has all fibres
isomorphic to $\mathfrak{A}_{2n}$. Notice that the induced vector
bundle $E':=U(E)\times_{U(n)}\Re^{2n}$ has an induced almost
symplectic structure and an induced almost complex structure.
Clearly $E'$ equipped with these structures is isomorphic to $(E,
\Omega, J)$. For each $\hbar\in (0,1]$ the $U(n)$-equivariant
linear map $\varphi_{\hbar}:\Re^{2n}\rightarrow \mathfrak{A}_{2n}$
also induces a bundle map $E\cong E'\rightarrow D$, which we still
denote by $\varphi_{\hbar}$.

\begin{definition} \label{quantization bundle:def}
Let $(E, \Pi, J)$ be as in Lemma~\ref{local unitary basis:lemma}.
We call the $C^*$-algebra bundle $D$ constructed above the
\emph{quantization bundle} of $(E, \Omega, J)$, and call the
bundle maps $\varphi_{\hbar}:E\rightarrow D$ the
\emph{quantization maps}.
\end{definition}

\begin{conventions} \label{bundle map:conventions}
Let $E$ and $D$ be a real vector bundle and a $C^*$-algebra bundle
over $M$ respectively. Every $\Re$-linear bundle map
$\varphi:E\rightarrow D$ extends to a $\Ce$-linear bundle map from
the complexified bundle $E\otimes \Ce$ to $D$ sending $f+ig$ to
$\varphi(f)+i\varphi(g)$, where $f, g\in \Gamma(E)$. We'll denote
this extended map still by $\varphi$.
\end{conventions}

Since $\mathfrak{A}_{2n}$ is unital the bundle $D$ contains the
trivial bundle $M\times \Ce$ as a subbundle naturally. Thus
$\Gamma_{\infty}(D)$ contains $C_{\infty}(M)$ as a subalgebra. We
are ready to construct strict quantizations for $M$:

\begin{theorem} \label{sq for almost sym:thm}
Let $(M, \Pi)$ be a symplectic manifold. Let $J$ be a compatible
almost complex structure on $T^*M$, and let $D$ and
$\varphi_{\hbar}$ be the quantization bundle and maps of $(T^*M,
\Pi, J)$. Let $\mathcal{A}_{\hbar}=\Gamma_{\infty}(D)$ and
$\pi_{\hbar}(f)=f+\hbar^{\frac{1}{2}}\varphi_{\hbar}(df)$ for all
$0<\hbar\le 1$ and $f\in C^{\infty}_c(M)$. Also let
$\mathcal{A}_0=C_{\infty}(M)$, and let $\pi_0$ be the canonical
embedding of $C^{\infty}_c(M)$ into $C_{\infty}(M)$. Let
$\{\mathcal{A}_{\hbar}\}$ be the subfield of the trivial
continuous field of $C^*$-algebras $[0,1]\times
\Gamma_{\infty}(D)$ over $[0,1]$. Then $\{\mathcal{A}_{\hbar},
\pi_{\hbar}\}$ is a faithful flabby strict quantization of $M$
over $I=[0,1]$ with $A=C^{\infty}_c(M)$. If a Lie group $G$ has an
action on $M$ preserving $\Pi$ and $J$, then this strict
quantization is $G$-equivariant.
\end{theorem}
\begin{proof}
We verify the conditions in Definition~\ref{sq:def}. Condition (1)
follows from our choice of $(\mathcal{A}_0, \pi_0)$. Condition (2)
follows from Definition~\ref{Heisenberg:def}(2)(3). Condition (3)
follows from Definition~\ref{Heisenberg:def}(1) and our discussion
at the beginning of this section. Thus $\{A_{\hbar},
\pi_{\hbar}\}$ is a strict quantization of $M$. Since
$A=C^{\infty}_c(M)$ this strict quantization is flabby. The
faithfulness follows from Lemma~\ref{faithful:lemma} below. In
fact we need only $\varphi_{\hbar}(\Re^{2n}\otimes \Ce)\cap \Ce
1_{\mathfrak{A}_{2n}}=\{0\}$ here, but we shall need the full
result of Lemma~\ref{faithful:lemma} later in Corollary~\ref{sq
not sdq:cor}.
\begin{lemma} \label{faithful:lemma}
Let $\mathfrak{A}_{2n}$ be a Heisenberg algebra. Let
$V_{\hbar}=\varphi_{\hbar}(\Re^{2n}\otimes
\Ce)+(\varphi_{\hbar}(\Re^{2n}\otimes \Ce))^*$. Then
$V_{\hbar}\cap \Ce 1_{\mathfrak{A}_{2n}}=\{0\}$ for every
$0<\hbar\le 1$.
\end{lemma}
\begin{proof} It suffices to show that the only $U(n)$-fixed
element in $V_{\hbar}$ is $0$. Endow $U(n)$ with the normalized
Haar measure. Let $\sigma$ be the canonical map
$\mathfrak{A}_{2n}\rightarrow (\mathfrak{A}_{2n})^{U(n)}$ defined
by $\sigma(a)=\int_{U(n)}\rho_h(a)\, dh$. Similarly define
$\tau:\Re^{2n}\otimes \Ce\rightarrow (\Re^{2n}\otimes
\Ce)^{U(n)}$. Then $\varphi_{\hbar}\circ\tau=\sigma\circ
\varphi_{\hbar}$. Clearly $(\Re^{2n})^{U(n)}=\{0\}$, and hence
$\tau(\Re^{2n}\otimes \Ce)=\{0\}$. Consequently
$(V_{\hbar})^{U(n)}=\sigma(V_{\hbar})=\{0\}$.
\end{proof}
Finally, the assertions about $G$-action is clear.
\end{proof}

The proof of Theorem~\ref{sq for almost sym:thm} generalizes to
Theorem~\ref{sq for almost sym bundle:thm} immediately:

\begin{theorem} \label{sq for almost sym bundle:thm}
Let $(M, \Pi)$ be an almost Poisson manifold. Let $(E, \Pi_E,
J_E)$ be as in Lemma~\ref{local unitary basis:lemma}. Suppose that
$\psi:T^*M\rightarrow E$ is an almost Poisson map
(Definition~\ref{alm sym bundle:def}). Let $D$ and
$\varphi_{\hbar}$ be the quantization bundle and maps of $(E,
\Pi_E, J_E)$. Let $\mathcal{A}_{\hbar}=\Gamma_{\infty}(D)$ and
$\pi_{\hbar}(f)=f+\hbar^{\frac{1}{2}}(\varphi_{\hbar}\circ
\psi)(df)$ for all $0<\hbar\le 1$ and $f\in C^{\infty}_c(M)$. Also
let $\mathcal{A}_0$ and $\varphi_0$ be as in Theorem~\ref{sq for
almost sym:thm}. Let $\{A_{\hbar}\}$ be the subfield of the
trivial continuous field of $C^*$-algebras $[0,1]\times
\Gamma_{\infty}(D)$ over $[0,1]$. Then $\{A_{\hbar},
\pi_{\hbar}\}$ is a faithful flabby strict quantization of $M$
over $I=[0,1]$ with $A=C^{\infty}_c(M)$. If a Lie group $G$ has an
action on $M$ preserving $\Pi$ and has an action on $(E, \Pi_E,
J_E)$ such that the projection $E\rightarrow M$ and the map $\psi$
are $G$-equivariant, then this strict quantization is
$G$-equivariant.
\end{theorem}

\section{Strict Quantizations of Almost Poisson Manifolds}
\label{sq alm P mani:sec}

Throughout the rest of this paper $(M, \Pi)$ will be an almost
Poisson manifold unless stated otherwise. To construct strict
quantizations of $M$ , by Theorem~\ref{sq for almost sym
bundle:thm} it suffices to find $(E, \Pi_E)$ and
$\psi:T^*M\rightarrow E$ as in Theorem~\ref{sq for almost sym
bundle:thm}. In fact there is a canonical choice of such $(E,
\Pi_E)$ and $\psi$. This follows from Lemmas~\ref{canonical
symplectic:lemma} and \ref{canonical embedding:lemma}, for which
we omit the proof.

\begin{lemma} \label{canonical symplectic:lemma}
Let $E\rightarrow M$ be a vector bundle, and let $E^*$ be its dual
bundle. Then $\Omega_x((u_1,v^*_1), (u_2, v^*_2))=<u_1,
v^*_2>-<u_2, v^*_1>$ defines an almost symplectic structure on
$E\oplus E^*$, where $u_j\in E_x, v^*_j\in E^*_x$ and $<\cdot,
\cdot>$ is the canonical pairing between $E$ and $E^*$.
\end{lemma}

\begin{lemma} \label{canonical embedding:lemma}
Let $(E, \Pi_E)$ be an almost Poisson bundle over $M$, and let
$\sigma: E\rightarrow E^*$ be the induced bundle map defined by
$\sigma_x(u)=\Pi_{E,x}(\cdot, u)$. Let $E\oplus E^*$ be endowed
with the canonical almost symplectic structure defined in
Lemma~\ref{canonical symplectic:lemma}. Then the bundle map
$\psi:E\rightarrow E\oplus E^*$ defined by
$\psi_x(u)=(\frac{1}{\sqrt{2}}u, \frac{1}{\sqrt{2}}\sigma_x(u))$
is an almost Poisson map.
\end{lemma}

By Theorem~\ref{sq for almost sym bundle:thm} we have:

\begin{theorem} \label{sq for almost Poisson manifold 1:thm}
Let $\Omega_M$ be the canonical almost symplectic structure on
$T^*M\oplus TM$ defined in Lemma~\ref{canonical symplectic:lemma}
for $E=T^*M$. Let $\psi:T^*M\rightarrow T^*M\oplus TM$ be the
bundle map defined in Lemma~\ref{canonical embedding:lemma} for
$(T^*M, \Pi)$. Let $J$ be a compatible almost complex structure
for $(T^*M\oplus TM, \Omega_M)$, and let $D$ and $\varphi_{\hbar}$
be the quantization bundle and maps of $(T^*M\oplus TM, \Omega_M,
J)$. Let $\{\mathcal{A}_{\hbar}, \pi_{\hbar}\}$ be as in
Theorem~\ref{sq for almost sym bundle:thm}. Then $\{A_{\hbar},
\pi_{\hbar}\}$ is a faithful flabby strict quantization of $M$
over $I=[0,1]$ with $A=C^{\infty}_c(M)$. If a Lie group $G$ has a
smooth  action on $M$ preserving $\Pi$ and $J$, then this strict
quantization is $G$-equivariant.
\end{theorem}

We are ready to prove Theorem~\ref{sq for almost Poisson manifold 2:thm}:

\begin{proof}[Proof of Theorem~\ref{sq for almost Poisson manifold 2:thm}]
Given any Riemannian metric on $M$, the induced isomorphism $TM\rightarrow T^*M$
induces an inner product on $T^*M$, and hence induces an inner product on $T^*M\oplus TM$ by
requiring $T^*M$ and $TM$ to be perpendicular to each other.
Therefore $T^*M\oplus TM$ has a $G$-invariant inner
product. Notice that given an inner product on an almost
symplectic bundle $(E, \Omega)$ there is a canonical way to
construct a compatible almost complex structure $J$ on $E$
\cite{Cannas01}. Thus $T^*M\oplus TM$ has a $G$-invariant compatible almost complex structure.
Now the assertions follow from Theorem~\ref{sq for almost Poisson manifold 1:thm}.
\end{proof}

\begin{corollary} \label{comp Lie:cor}
For any smooth action of a compact Lie group $G$ on $M$ preserving
$\Pi$, there is a $G$-equivariant faithful strict quantization of
$M$ over $I=[0,1]$ with $A=C^{\infty}_c(M)$.
\end{corollary}
\begin{proof}
For any smooth action of a compact Lie group, the manifold
admits an invariant Riemannian metric  by integrating any given
Riemannian metric.
\end{proof}

Rieffel showed \cite{Rieffel89} that there is no strict
deformation quantization of the rotationally invariant symplectic
structure on $S^2$  respecting the action of $SO(3)$. So this
gives us some sign on how more restrictive strict deformation
quantizations are  than strict quantizations.

\begin{corollary} \label{dual:cor}
If a Lie group $G$ has a bi-invariant Riemannian metric,  then $\mathfrak{g}^*$
equipped with the Lie-Poisson bracket \cite{CW01}
admits a faithful strict quantization with
$A=C^{\infty}_c(\mathfrak{g}^*)$ equivariant under the coadjoint action of
$G$,  where $\mathfrak{g}$ is the Lie algebra of $G$ and
$\mathfrak{g}^*$ is the dual.
\end{corollary}
\begin{proof}
Identify the cotangent space of $\mathfrak{g}^*$ at each point
with $\mathfrak{g}$.  Then for any $\xi\in \mathfrak{g}^*$ and
$g\in G$,  the isomorphism $T^*_{\xi}\mathfrak{g}^*\to
T^*_{Ad^*_g(\xi)}\mathfrak{g}^*$ is exactly $Ad_g:\mathfrak{g}\to
\mathfrak{g}$.  In particular considering $\xi=0$, we see that
$\mathfrak{g}^*$ admits an invariant Riemannian metric if and only
if the vector space $\mathfrak{g}$ has an inner product invariant
under the adjoint action of $G$, if and only if $G$ has a
bi-invariant Riemannian metric.
\end{proof}

\begin{example} \label{trivial TM:eg}
Assume that $TM$ is trivial. Let $X_1, \cdots, X_m\in \Gamma(TM)$ giving the trivialization of
$TM$. Let $X^*_1, \cdots, X^*_m\in \Gamma(T^*M)$ be the dual basis.
Define $J$ on $T^*M\oplus TM$ by $J(X^*_k)=X_k$ and $J(X_k)=-X^*_k$ for
$1\le k\le m$. Then the quantization bundle
$D$ is the trivial bundle $M\times \mathfrak{A}_{2m}$, and thus
$\Gamma_{\infty}(D)=C_{\infty}(M, \mathfrak{A}_{2m})=
C_{\infty}(M)\otimes \mathfrak{A}_{2m}$. Let
\begin{eqnarray*}
\beta_{jk}=\Pi(X^*_j, X^*_k)
\end{eqnarray*}
be the \emph{structure functions}. Then the bundle map $\sigma:T^*M\rightarrow TM$ is determined by
$\sigma(X^*_k)=\sum^m_{j=1}\beta_{jk}X_j$.
Thus for any $f, g\in C^{\infty}_c(M)$ we have
\begin{eqnarray*}
\{f,g\}&=&\sum_{1\le j,k\le m}\beta_{jk}X_j(f)X_k(g),\\
\pi_{\hbar}(f)&=&f\otimes 1+\frac{1}{\sqrt{2}}\hbar^{\frac{1}{2}}(\sum_{1\le j\le m}X_j(f)\otimes \varphi_{\hbar}(e_j)+
\sum_{1\le j,k\le m}\beta_{jk}X_k(f)\otimes \varphi_{\hbar}(e_{m+j})),
\end{eqnarray*}
where $e_j$ and $\varphi_{\hbar}$ are as in Definition~\ref{Heisenberg:def}.
\end{example}

\begin{example} \label{tori:eg}
Let $M$ be the $m$-torus $\mathbb{T}^m$, and let $x_1, \cdots,
x_m$ be the standard coordinates. Let $\theta$ be a real
skew-symmetric $m\times m$ matrix. Define a Poisson bracket
$\{\cdot,\cdot\}$ on $M$ by
\begin{eqnarray*}
\{f, g\}=\frac{1}{2}\sum_{1\le j, k\le m}\theta_{jk}\frac{\partial f}{\partial
x_j}\frac{\partial g}{\partial x_k}.
\end{eqnarray*}
Let $X_j=\frac{\partial }{\partial x_j}$.
Let $\beta_{jk}$ and $J$ be as in Example~\ref{trivial TM:eg}.
Then $\beta_{jk}=\frac{1}{2}\theta_{jk}$. Thus
\begin{eqnarray*}
\pi_{\hbar}(f)=f\otimes 1+
\frac{1}{\sqrt{2}}\hbar^{\frac{1}{2}}(\sum_{1\le j\le m}\frac{\partial f}{\partial x_j}\otimes \varphi_{\hbar}(e_j)+
\frac{1}{2}\sum_{1\le j,k\le m}\theta_{jk}\frac{\partial f}{\partial x_k}\otimes \varphi_{\hbar}(e_{m+j})).
\end{eqnarray*}
This is very different from Rieffel's Moyal product approach \cite{Rieffel89, Rieffel93},
which leads to the noncommutative torus
$A_{\theta}$.
\end{example}

\begin{example} \label{dual:eg}
Let $\mathfrak{g}$ be a Lie algebra, and let $M$ be the dual $\mathfrak{g}^*$ equipped with the Lie-Poisson bracket.
Let $v_1, \cdots, v_m$ be a basis of $\mathfrak{g}$, and let
$\mu_1, \cdots, \mu_m$ be the dual basis of $\mathfrak{g}^*$.
Let $c_{jkl}$ be the \emph{structure constants} satisfying $[v_j, v_k]=\sum c_{jkl}v_l$.
We may take $X_j$ in Example~\ref{trivial TM:eg} to be $\mu_j$.
Then $X^*_j=v_j$, and $\beta_{jk}=\sum_lc_{jkl}v_l$. Thus
\begin{eqnarray*}
\pi_{\hbar}(f)=f\otimes 1+\frac{1}{\sqrt{2}}\hbar^{\frac{1}{2}}(\sum_{1\le j\le m}\mu_j(f)\otimes \varphi_{\hbar}(e_j)+
\sum_{1\le j,k,l\le m}c_{jkl}v_l\mu_k(f)\otimes \varphi_{\hbar}(e_{m+j})).
\end{eqnarray*}
\end{example}
\section{Homotopy}
\label{homotopy:sec}

Our construction in Theorem~\ref{sq for almost Poisson manifold
1:thm} depends on the choice of $J$. We define homotopy of strict
quantizations first, then show that the homotopy class of our
construction in independent of the choice of $J$
(Proposition~\ref{homotopy:prop}). The definition of homotopy of
strict quantizations is similar to the usual definition of
homotopy of homomorphisms between $C^*$-algebras:

\begin{definition} \label{homotopy:def}
Let $\{\mathcal{A}^j_{\hbar}, \pi^j_{\hbar}\}$ be strict
quantizations of $(M, \Pi)$ over $I$ for $A$, where $j=0, 1$. By a
\emph{homotopy} of these two strict quantizations, we mean a
continuous field of $C^*$-algebras $\{\mathcal{A}_{\hbar, t}\}$
over $I\times [0,1]$ and linear maps $\pi_{\hbar,t}:A\rightarrow
\mathcal{A}_{\hbar, t}$ such that

(1) the restriction of this field on $I\times \{t\}$ gives a
strict quantization of $(M, \Pi)$ over $A$ for each $t\in [0,1]$
in a uniform way, \ie for all $f,g\in A$ we have
\begin{eqnarray*}
   \lim_{\hbar\to 0}\sup_{0\le t\le 1}\pa
   [\pi_{\hbar,t}(f),\pi_{\hbar,t}(g)]/(i\hbar)-\pi_{\hbar,t}(\{f,g\})\pa
   = 0,
\end{eqnarray*}

(2) for $t=0,1$ the restriction of this field gives the strict
quantizations  $\{\mathcal{A}^0_{\hbar}, \pi^0_{\hbar}\}$ and
$\{\mathcal{A}^1_{\hbar}, \pi^1_{\hbar}\}$ respectively.
\end{definition}

\begin{remark} \label{homotopy:remark}
One may also define a weaker notion of homotopy without requiring
the convergence in (1) to be uniformly. We adopt the stronger one
because the homotopies we construct here all satisfy the uniform
condition.
\end{remark}

Clearly homotopy is an equivalence relation between strict
quantizations of $M$.

\begin{proposition} \label{homotopy:prop}
The homotopy class of the strict quantization in Theorem~\ref{sq
for almost Poisson manifold 1:thm} does not depend on the choice
of the compatible almost complex structure $J$ on $T^*M\oplus TM$.
\end{proposition}
\begin{proof}
Let $J_0$ and $J_1$ be two compatible almost complex structures on
$T^M\oplus TM$. Let $<\cdot, \cdot>_0$ and $<\cdot, \cdot>_1$ be
the induced inner products. Let $<\cdot, \cdot>_t=t<\cdot,
\cdot>_1+(1-t)<\cdot, \cdot>_0$ for $0\le t\le 1$. The canonical
way of constructing compatible almost complex structure from given
inner product \cite{Cannas01} is continuous. Thus we get a
continuous family of compatible almost complex structures $J_t$ on
$T^*M\oplus TM$. Let $(D_t, \varphi_{\hbar,t})$ be the
corresponding quantization bundle and maps. Denote
$\Gamma_c(T^*M\oplus TM)$ the space of compactly supported
sections of $T^M\oplus TM$. Then the sections
$(\varphi_{\hbar,t}(f))_{0\le t\le 1}$ for $0<\hbar\le 1$ and
$f\in \Gamma_c(T^*M\oplus TM)$ generate a continuous field of
$C^*$-algebras over $[0,1]$ with fibre $\Gamma_{\infty}(D_t)$ at
$t$. Now it is easy to see that the strict quantizations in
Theorem~\ref{sq for almost Poisson manifold 1:thm} associated with
$J_t$'s combine together to give a homotopy of the ones associated
with $J_0$ and $J_1$.
\end{proof}

For given $(M, \Pi, A, I)$ we don't know whether there is only one homotopy
class of strict quantizations. But this is the case when
$\Pi=0$ and $I$ is an interval. Notice first that for $\Pi=0$
there is a canonical trivial strict quantization over any $I$,
namely $\mathcal{A}_{\hbar}=C_{\infty}(M)$ and $\pi_{\hbar}$ is the
canonical inclusion of $A$ into $C_{\infty}(M)$.
\begin{proposition} \label{Pi=0:prop}
When $\Pi=0$ and $I$ is an interval, every strict quantization
over $I$ is homotopic to the canonical one.
\end{proposition}
\begin{proof}
Let $\{\mathcal{A}_{\hbar}, \varphi_{\hbar}\}$ be a strict
quantization. Define a map $\gamma:I\times [0,1]\rightarrow I$ by
$\gamma(\hbar, t)=t\hbar$. Then the pull back of the field
$\{\mathcal{A}_{\hbar}\}$ under $\gamma$ is a continuous field
over $I\times [0,1]$ with fibre $\mathcal{A}_{\hbar,
t}=\mathcal{A}_{t\hbar}$ at $(\hbar, t)$. For each $f\in A$ let
$\{\pi_{\hbar, t}(f)\}$ be the pull back of the section
$\{\pi_{\hbar}(f)\}$ under $\gamma$, namely
$\pi_{\hbar,t}(f)=\pi_{t\hbar}(f)$. Then clearly this is a
homotopy between the canonical strict quantization and
$\{\mathcal{A}_{\hbar}, \varphi_{\hbar}\}$.
\end{proof}

\section{Local Strict Quantizations}
\label{local strict quantization:sec}

There are two different meanings for a strict quantization to be
local. The first one is an intuitive one, meaning that the
construction is local in the sense that we construct strict
quantizations for open subsets of $M$ first, then gluing them
together to get  strict quantization for $M$. This includes our
construction in Theorem~\ref{sq for almost sym bundle:thm} and the
construction in \cite{NNP02}. The second one means that the
algebras and maps $\{\mathcal{A}_{\hbar}, \pi_{\hbar}\}$ are local
in the sense that $\mathcal{A}_{\hbar}\subseteq
\mathcal{D}_{\hbar}$ for some (upper-semi)continuous field of
$C^*$-algebras $D_{\hbar}$ and the maps
$\pi_{\hbar}:C_{\infty}(M)=\Gamma_{\infty}(M\times \Ce)\rightarrow
\mathcal{A}_{\hbar}\hookrightarrow \Gamma_{\infty}(D_{\hbar})$ are
fibrewise. Here we'll concentrate on the second meaning.

Let $X$ be a locally compact Hausdorff space. Recall that an {\it
$C_{\infty}(X)$-algebra} is a $C^*$-algebra $\mathcal{A}$ with an
injective nondegenerate homomorphism
$\gamma:C_{\infty}(X)\rightarrow \mathcal{M}(\mathcal{A})$ such
that $\gamma(C_{\infty}(X))$ being contained in the center
$\mathcal{Z}\mathcal{M}(\mathcal{A})$ of the multiplier algebra
$\mathcal{M}(\mathcal{A})$ \cite{Kasparov88}. This is equivalent
to saying that $\mathcal{A}$ is the global section algebra of an
upper-semicontinuous field of $C^*$-algebras over $X$
\cite{Nilsen98}. Under this correspondence the fibre algebra of
the field at $x\in X$ is
$\mathcal{A}/\overline{\gamma(I_x)\mathcal{A}}$, where $I_x=\{h\in
C_{\infty}(X):h(x)=0\}$. This motivates our definition of local
strict quantizations:

\begin{definition} \label{local sq:def}
Let $\{\mathcal{A}_{\hbar}, \pi_{\hbar}\}$ be a strict
quantization of $(M, \Pi)$ on $I$. We call $\{\mathcal{A}_{\hbar},
\pi_{\hbar}\}$ \emph{local} if each $\mathcal{A}_{\hbar}$ is a
$C_{\infty}(M)$-algebra with
$\gamma_{\hbar}:C_{\infty}(M)\hookrightarrow
\mathcal{Z}\mathcal{M}(\mathcal{A}_{\hbar})$ such that $\pa
\pi_{\hbar}(f)-\gamma_{\hbar}(f)\pa\rightarrow 0$ as
$\hbar\rightarrow 0$ for every $f\in A$.
\end{definition}

Clearly the strict quantizations in Theorem~\ref{sq for almost sym bundle:thm} are local.
But the ones in \cite{NNP02} are not.

\begin{proposition} \label{local sq:prop}
Let $\{\mathcal{A}_{\hbar}, \pi_{\hbar}\}$ be a local strict quantization of
$(M, \Pi)$ on $I$. Let $f, g\in A_{sa}$. If $\pi_{\hbar}(f),
\pi_{\hbar}(g)\in (\mathcal{A}_{\hbar})_{sa}$ for all $\hbar\in I$, then
$\{f, g\}=0$.
\end{proposition}
\begin{proof} Using the embeddings $\gamma_{\hbar}$ in
Definition~\ref{local sq:def} we'll identify $C_{\infty}(M)$ as a
subalgebra of $\mathcal{M}(\mathcal{A}_{\hbar})$. Then $\pa
\pi_{\hbar}(\{f,g\})-\{f,g\}\pa \rightarrow 0$ as
$\hbar\rightarrow 0$. Thus Definition~\ref{sq:def}(3) becomes
\begin{eqnarray} \label{(f,g):eq}
   \lim_{\hbar\to 0}\pa
   [\pi_{\hbar}(f), \pi_{\hbar}(g)]/(i\hbar)-\{f,g\}\pa
   = 0.
\end{eqnarray}
For $x\in M$ let $I_x=\{h\in C_{\infty}(M):h(x)=0\}$, and let
$\mathcal{A}_{\hbar,x}=\mathcal{A}_{\hbar}/\overline{I_x\mathcal{A}_{\hbar}}$.
Let $\beta_{\hbar,x}:\mathcal{A}_{\hbar}\rightarrow
\mathcal{A}_{\hbar,x}$ be the quotient map. Notice that the
identity of $\mathcal{M}(\mathcal{A}_{\hbar, x})$ is
$\beta_{\hbar,x}(h)$ for any $h\in C_{\infty}(M)$ with $h(x)=1$.
Taking $\beta_{\hbar,x}$ on (\ref{(f,g):eq}) we get
\begin{eqnarray} \label{(f,g)2:eq}
   \lim_{\hbar\to 0}\pa
   [\beta_{\hbar,x}(\pi_{\hbar}(f)),\beta_{\hbar,x}(\pi_{\hbar}(g))]/(i\hbar)-\{f,g\}(x)\pa
   = 0.
\end{eqnarray}
Notice that
$[\beta_{\hbar,x}(\pi_{\hbar}(f)),\beta_{\hbar,x}(\pi_{\hbar}(g))]/(i\hbar)$
is a {\it self-commutator}, \ie of the form $[S^*,S]$ for some $S$
(for instance
$S=(2\hbar)^{-1/2}\beta_{\hbar,x}(\pi_{\hbar}(f))-i(2\hbar)^{-1/2}\beta_{\hbar,x}(\pi_{\hbar}(g))$).
It is known that self-commutators can't be invertible
\cite[Corollary 1]{Radjavi66}. Thus $\{f,g\}(x)=0$.
\end{proof}

\begin{corollary} \label{sa local sq:cor}
An almost Poisson manifold $(M, \Pi)$ admits a Hermitian local
strict quantization if and only if $\Pi=0$.
\end{corollary}
\begin{proof} Assume that $\Pi\neq 0$ and that $M$ admits a
Hermitian local strict quantization. Then we can find a covector
field $Y^*\in \Gamma(T^*M)$ with $\sigma(Y^*)\neq 0$, where
$\sigma$ is as in Lemma~\ref{canonical embedding:lemma} for
$E=T^*M$. We claim that there is a vector field $X\neq 0$ such
that $X(df)=0$ for all $f\in A$. If $\sigma(dg)\neq 0$ for some
$g\in A_{sa}$, by Proposition~\ref{local sq:prop} we may take
$X=\sigma(dg)$. Otherwise we may take $X=\sigma(Y^*)$. Let $Z$ be
a nonconstant integral curve of $X$. Then the restriction of every
$f\in A$ on $Z$ is constant, which contradicts $A$ being dense in
$C_{\infty}(M)$. This proves the "only if" part. The "if" part is
trivial.
\end{proof}

\begin{remark} \label{sa local sq:remark}
We don't know when a local strict quantization is homotopic to a Hermitian strict quantization.
As a comparison, a star product on a symplectic manifold is equivalent to a Hermitian one
if and only if its characteristic class is Hermitian \cite{Neumaier99}.
Also every Poisson manifold has Hermitian star products \cite{BW01}.
\end{remark}

Thus the strict quantizations in Theorem~\ref{sq for almost sym bundle:thm} are not Hermitian unless
$\Pi=0$. In fact we can say more:

\begin{proposition} \label{sq not sdq:prop}
Let $\{\mathcal{A}_{\hbar}, \pi_{\hbar}\}$ be a local strict
quantization of $(M, \Pi)$ on $I$. Identify $C_{\infty}(M)$ as a
subalgebra of $\mathcal{M}(\mathcal{A}_{\hbar})$ via
$\gamma_{\hbar}$. Let $\tau_{\hbar}(f)=\pi_{\hbar}(f)-f$ for
$0<\hbar\le 1$ and $f\in A$. Assume that
\begin{eqnarray} \label{sq not sdq:eq}
A\cap (\tau_{\hbar}(A)+(\tau_{\hbar}(A))^*)=\{0\}
\end{eqnarray}
for every $0<\hbar\le 1$. If this is a strict deformation
quantization, then it is Hermitian and $\Pi=0$.
\end{proposition}
\begin{proof}
Let $f\in A$ and $0<\hbar\le 1$. Then
$(\pi_{\hbar}(f))^*=\pi_{\hbar}(g)$ for some $g\in A$. Thus
$f^*+(\tau_{\hbar}(f))^*=g+\tau_{\hbar}(g)$. By our assumption
$f^*=g$. Thus this strict quantization is Hermitian. By Corollary~\ref{sa
local sq:cor} $\Pi=0$.
\end{proof}

\begin{corollary} \label{sq not sdq:cor}
Let $\{\mathcal{A}_{\hbar}, \pi_{\hbar}\}$ be a local strict
quantization of $(M, \Pi)$ on $I$ with $\Pi\neq 0$. Let
$\tau_{\hbar}$ be as in Proposition~\ref{sq not sdq:prop} and
assume that (\ref{sq not sdq:eq}) holds for all $0<\hbar\le 1$.
Then $\{\mathcal{A}_{\hbar}, \pi_{\hbar}\}$ can't be restricted to
a dense $*$-subalgebra of $A$ to get a strict deformation
quantization of $(M, \Pi)$. In particular, the strict
quantizations in Theorem~\ref{sq for almost sym bundle:thm} can't
be restricted to a dense $*$-subalgebra of $C^{\infty}_c(M)$ to
get a strict deformation quantization unless $\Pi=0$.
\end{corollary}
\begin{proof}
By Lemma~\ref{faithful:lemma} the strict quantizations in
Theorem~\ref{sq for almost sym bundle:thm} satisfy (\ref{sq not
sdq:eq}).
\end{proof}

\begin{corollary} \label{sq not sdq 2:cor}
For any $(M, \Pi)$ there is a faithful flabby strict
quantization, which can't be restricted to  any dense
$*$-subalgebra of $A$ to get a strict deformation quantization.
\end{corollary}
\begin{proof} The case $\Pi\neq 0$ follows from Theorem~\ref{sq for almost Poisson manifold 1:thm}
and Corollary~\ref{sq not sdq:cor}. The case $\Pi=0$ is settled in
\cite{Li7}.
\end{proof}

Corollary~\ref{sq not sdq 2:cor} gives Question 25 in \cite{Rieffel98a} a negative answer,  which
asks whether there is an example of a faithful strict quantization
of a Poisson manifold such that it's impossible to restrict it to some dense
$*$-subalgebra to get a strict deformation quantization.
This leaves the question whether we can require the strict quantization
to be Hermitian.

\section{Functorial Properties}
\label{functor:sec}

It is unlikely that there is a universal way to construct a canonical
strict quantization for each Poisson manifolds such that
it gives a contravariant functor from the category of Poisson manifolds
with (proper) Poisson maps to the category of continuous field of $C^*$-algebras over
$I$ \cite{Landsman01}. Instead, Landsman proposed other categories closely related
to Morita equivalence, and showed that there is such a functor on the subcategory of
dual of integrable Lie algebroids.

Though our construction in Theorem~\ref{sq for almost Poisson manifold 1:thm} doesn't
give a contravariant functor, it does has some properties similar to functors.
In this section we discuss two questions:

(1) fixing a strict quantization $\{\mathcal{A}^M_{\hbar},
\pi^M_{\hbar}\}$ of $(M, \Pi_M)$ on $I$, for any proper Poisson
map $\phi$ from $M$ to another almost Poisson manifold $(N,
\Pi_N)$ can we find a strict quantization
$\{\mathcal{A}^N_{\hbar}, \pi^N_{\hbar}\}$ of $(N, \Pi_N)$ on $I$
with a "homomorphism" of these two strict quantizations extending
$\phi^*:C_{\infty}(N)\rightarrow C_{\infty}(M)$, \ie a
homomorphism $\xi_{\hbar}:\mathcal{A}^N_{\hbar}\rightarrow
\mathcal{A}^M_{\hbar}$ for each $\hbar\in I$ such that these maps
$\{\phi_{\hbar}\}$ send continuous sections to continuous ones and
$\xi_{\hbar}\circ \pi^N_{\hbar}=\pi^M_{\hbar}\circ \phi^*$?

(2) the similar question but fixing the strict quantization of $N$ instead.

The first question has a positive answer because of Theorem~\ref{sq
for almost Poisson manifold 2:thm} and the following Proposition,
whose proof is just routine verification.

\begin{proposition} \label{functor fixing M:prop}
Let $\{\mathcal{A}^M_{\hbar}, \pi^M_{\hbar}\}$  and
$\{\mathcal{A}^N_{\hbar}, \pi^N_{\hbar}\}$ be strict quantizations
of $(M, \Pi_M)$  and $(N, \Pi_N)$ on $I$ for $A^M$ and $A^N$
respectively. Let $\phi:M\rightarrow N$ be a proper Poisson map
with $\phi^*(A^N)\subseteq A^M$. Then the sections $\{(f, g): f\in
\Gamma(\{\mathcal{A}^M_{\hbar}\}), g\in
\Gamma(\{\mathcal{A}^N_{\hbar}\}), f_0=\phi^*(g_0)\}$ determine a
continuous field of $C^*$-algebras over $I$ with fibre
$\mathcal{A}^M_{\hbar}\oplus \mathcal{A}^N_{\hbar}$ at $\hbar\neq
0$ and fibre $C_{\infty}(M)$ at $\hbar=0$. And
$(\pi^M_{\hbar}\circ \phi^*)\oplus \pi^N_{\hbar}$ for $\hbar\neq
0$ give a strict quantization for $(N, \Pi_N)$ on $I$ for $A^N$.
\end{proposition}

\begin{remark} \label{functor fixing M:remark}
(1) When both $\{\mathcal{A}^M_{\hbar}, \pi^M_{\hbar}\}$  and
$\{\mathcal{A}^N_{\hbar}, \pi^N_{\hbar}\}$ are local in the sense
of Definition~\ref{local sq:def}, so is
$\{\mathcal{A}^M_{\hbar}\oplus \mathcal{A}^N_{\hbar},
(\pi^M_{\hbar}\circ \phi^*)\oplus \pi^N_{\hbar}\}$;

(2) When $(M, \Pi_M)=(N, \Pi_N)$ and $\phi=id_M$, Proposition~\ref{functor fixing
M:prop} shows that the set of isomorphism classes of strict
quantizations of $(M, \Pi)$ over $I$ for $A$ has a natural abelian
semigroup structure. Clearly the addition is compatible with homotopy
defined in Definition~\ref{homotopy:def}. Thus the set of homotopy
classes of strict quantizations of $(M, \Pi)$ over $I$ for $A$ is
also an abelian semigroup.
\end{remark}

For the second question we have a partial positive answer:

\begin{proposition} \label{functor fixing N:prop}
Let $(N, \Pi_N)$ be an almost symplectic manifold, and let $k\ge \dim N$.
Then there is a strict quantization $\{\mathcal{A}^N_{\hbar},
\pi^N_{\hbar}\}$ of $(N, \Pi_N)$ as constructed in Theorem~\ref{sq for almost Poisson manifold 1:thm}
such that
for any proper Poisson map $\phi:(M, \Pi_M)\rightarrow (N, \Pi_N)$ with $k\ge \dim M$ there is a
strict quantization $\{\mathcal{A}^M_{\hbar}, \pi^M_{\hbar}\}$ of $(M, \Pi_M)$
as constructed in Theorem~\ref{sq for almost Poisson manifold 1:thm}
and homomorphisms $\xi_{\hbar}:\mathcal{A}^N_{\hbar}\rightarrow \mathcal{A}^M_{\hbar}$ sending
continuous sections to continuous ones with $\xi_{\hbar}\circ \pi^N_{\hbar}=\pi^M_{\hbar}\circ \phi^*$.
\end{proposition}
\begin{proof} Let $n=\dim N$, and let $n\le m\le k$.
We'll choose a special asymptotic representation of Heisenberg
commutation relations of dimension $2m$. Let $e_1, \dots, e_{2m}$
and $e'_1, \cdots, e'_{2k}$ be the standard basis of $\Re^{2m}$
and $\Re^{2k}$ respectively. Then the linear map $\eta:
\Re^{2m}\rightarrow \Re^{2k}$ defined by $\eta(e_j)=e'_j,\,
\eta(e_{j+m})=e'_{j+k}$ for $1\le j\le m$ preserves the standard
symplectic structure and the standard compatible almost complex
structure. Thus $U(m)$ can be thought of as the subgroup of $U(k)$
fixing $e'_l, e'_{l+k}$ for $m<l\le k$. Let $(\mathfrak{A}_{2k},
\varphi_{\hbar})$ be an asymptotic representation of Heisenberg
commutation relations of dimension $2k$. Let $\mathfrak{A}_{2m}$
be the $C^*$-subalgebra generated by $\cup_{0<\hbar\le
1}(\varphi_{\hbar}\circ \eta)(\Re^{2m})$. Then
$(\mathfrak{A}_{2m}, \varphi_{\hbar}\circ \eta)$ is an asymptotic
representation of Heisenberg commutation relations of dimension
$2m$.

Fix a compatible almost complex structure $J^N$ on $T^*N\oplus
TN$. Let $x\in M$. Since $(N, \Pi_N)$ is almost symplectic,
$(T^*\phi)_{\phi(x)}$ is injective. Thus $n\le \dim M$. Let $(D^N,
\sigma^N, \varphi^N_{\hbar})$ and $\sigma^M$ be as in
Theorem~\ref{sq for almost Poisson manifold 1:thm} and
Lemma~\ref{canonical embedding:lemma} for $N$ and $M$
respectively. Since $(N, \Pi_N)$ is almost symplectic, $\sigma^N$
is invertible. Then we have a linear map
$\theta_x:=\sigma^M_x\circ (T^*\phi)_{\phi(x)}\circ
(\sigma^N_{\phi(x)})^{-1}:TN_{\phi(x)}\rightarrow TM_x$. Let
$\zeta_x:=(T^*\phi)_{\phi(x)}\oplus \theta_x: T^*N_{\phi(x)}\oplus
TN_{\phi(x)}\rightarrow  T^*M_x\oplus TM_x$. Easy computation
shows that $TM_x\circ \theta_x$ is the identity map on
$TN_{\phi(x)}$, and that $\zeta_x$ preserves the canonical
symplectic structure on $T^*N_x\oplus TN_x$. Then
$\{\zeta_x(T^*N_{\phi(x)}\oplus TN_{\phi(x)})\}$ is an almost
symplectic subbundle of $T^*M\oplus TM$, which we'll denote by
$E$. Let $F_x=(E_x)^{\perp}$ with respect to the almost symplectic
structure. Then $\{F_x\}$ is also an almost symplectic subbundle
of $T^*M\oplus TM$, which we'll denote by $F$. Clearly $T^*M\oplus
TM=E\oplus F$, and $J^N$ induces a compatible almost complex
structure $J^E$ on $E$. Take a compatible almost complex structure
$J^F$ on $F$. Then $J^M:=J^E\oplus J^F$ is a compatible almost
complex structure on $T^*M\oplus TM$. Let $(D^M,
\varphi^N_{\hbar})$ be as in Theorem~\ref{sq for almost Poisson
manifold 1:thm} for $(M, \Pi_M, J^M)$.

Say $m=\dim M$.
Notice that
\begin{eqnarray*}
\zeta_x(u_1), \cdots, \zeta_x(u_n), \mu_1, \cdots, \mu_{m-n},
\zeta_x(v_1), \dots, \zeta_x(v_n), \gamma_1, \cdots, \gamma_{m-n}
\end{eqnarray*}
is a unitary basis of $T^*M_x\oplus TM_x$ for any unitary basis
$u_1, \cdots, u_n, v_1, \cdots, v_n$ of $T^*N_{\phi(x)}\oplus
TN_{\phi(x)}$ and any unitary basis $\mu_1, \cdots,
\mu_{m-n},\gamma_1, \cdots, \gamma_{m-n}$ of $F_x$. Because of our
choice of $\mathfrak{A}_{2m}$ and $\mathfrak{A}_{2n}$ the map
$\zeta_x$ determines a unital $C^*$-algebra embedding
$\xi_x:D^N_{\phi(x)}\rightarrow D^M_x$ such that
$\varphi^M_{\hbar,x}\circ \zeta_x=\xi_x \circ
\varphi^N_{\hbar,\phi(x)}$ for all $0<\hbar\le 1$. Since $\phi$ is
proper, the $\xi_x$'s combine to give a homomorphism
$\xi:\Gamma_{\infty}(D^N)\rightarrow \Gamma_{\infty}(D^M)$ whose
restriction on $C_{\infty}(N)$ is $\psi^*$. Clearly
$\xi_{\hbar}:=\xi$ satisfy the requirement.
\end{proof}

\begin{remark} \label{functor fixing N:remark}
(1) When
both $M$ and $N$ have $G$-equivariant Riemannian metrics
there is an obvious $G$-equivariant version of Proposition~\ref{functor fixing N:prop};

(2) If we can find a $C^*$-algebra $\mathfrak{A}_{\infty}$ with
an action of $U(\infty):=\cup_{n\in \Ne}U(n)$ and linear maps
$\varphi^n_{\hbar}:\Re^{2n}\rightarrow \mathfrak{A}_{\infty}$ for all $n\in \Ne$
compatible with the embedding $\eta:\Re^{2n}\rightarrow \Re^{2k}$ in the proof of Proposition~\ref{functor fixing N:prop}
such that for each $n$ these maps give an asymptotic representation of Heisenberg commutation relations, then
we can use $\mathfrak{A}_{\infty}$ instead of $\mathfrak{A}_{2k}$ in the above proof and hence throw away the requirement
$k\ge \dim M$. But we don't know whether such infinite dimensional asymptotic representation of Heisenberg commutation relations
exists or not.
\end{remark}

\section{Asymptotic Representation of Heisenberg Commutation Relations}
\label{asym rep:sec}

In this section we prove Theorem~\ref{asymptotic rep exit:thm}.

\begin{lemma} \label{asym rep 1:lemma}
Let $H$ be a separable Hilbert space with orthonormal basis
$\{e_j\}^{\infty}_{j=1}\cup\{e'_j\}^{\infty}_{j=1}$. Then there
exist norm continuous paths $T(\hbar), S(\hbar)$ of operators in
$B(H)$ for $0< \hbar\le 1$ such that

(1) $[T(\hbar), S(\hbar)](e_{2j-1})=(1+\hbar)ie_{2j-1}$,

(2) $[T(\hbar), S(\hbar)](e_{2j})=(1-\hbar)ie_{2j}$,

(3) $[T(\hbar), S(\hbar)](e'_j)=ie'_j$,

(4) $\hbar\pa T(\hbar)\pa$ and $\hbar\pa S(\hbar)\pa$ are
bounded uniformly in $\hbar$.

\end{lemma}

Brown and Pearcy \cite{BP65} proved that for a separable Hilbert
space $H$ an operator $R\in B(H)$ is a commutator if and only if
it is not a non-zero  scalar modulo compact operators. For the
"if" part, their proof is constructive. Since we need $T(\hbar)$
and $S(\hbar)$ to be continuous, and want some control on their
norms, and the construction in \cite{BP65} depends on some choices
of isomorphisms of Hilbert spaces, we write down the proof of
Lemma~\ref{asym rep 1:lemma} here, though it is just following the
construction in \cite{BP65}.

\begin{proof}
Let $\eta_j=\frac{1}{\sqrt{2}}(e_{2j-1}+e_{2j}),
\eta'_j=\frac{1}{\sqrt{2}}(e_{2j-1}-e_{2j})$.  Define $R(\hbar)\in
B(H)$ by $R(\hbar)(e_{2j-1})=(1+\hbar)ie_{2j-1},\,
R(\hbar)(e_{2j})=(1-\hbar)ie_{2j}$ and $R(\hbar)(e'_j)=ie'_j$.
Define $Z(\hbar)\in B(H)$ by $Z(\hbar)(\eta_j)=\eta_j,\,
Z(\hbar)(e'_i)=e'_j$ and $Z(\hbar)(\eta'_j)=\eta_j+\hbar\eta'_j$.
Then on $span_{\Ce}\{e_{2j-1}, e_{2j}\}$, with $\eta_j,\eta'_j$
as basis,  $Z(\hbar)$ and $Z(\hbar)^{-1}$ are $\left(
\begin{array}{cr} 1 & 1
\\ 0 &  \hbar \end{array} \right)$ and $\hbar^{-1}\left(\begin{array}{cr}  \hbar & -1 \\
0 & 1\end{array} \right)$ respectively. Therefore $\pa Z(\hbar)\pa
< 3$ and $\pa Z(\hbar)^{-1}\pa < 3\hbar^{-1}$. Hence it suffices
to find continuous paths $T(\hbar), S(\hbar)$ such $[T(\hbar),
S(\hbar)]=Z(\hbar)^{-1}R(\hbar)Z(\hbar)$  and $\pa T(\hbar)\pa,\,
\pa S(\hbar)\pa $ are bounded uniformly in $\hbar$.

If we identify $H$ with the closure span of
$\{e'_i\}^{\infty}_{1}, \{\eta'_i\}^{\infty}_{1},
\{\eta_i\}^{\infty}_{1}$ respectively  and hence identify $H$ with
$H\oplus H\oplus H$, simple calculation shows that
\begin{eqnarray*}
Z(\hbar)^{-1}R(\hbar)Z(\hbar)=i\begin{pmatrix} I & 0 & 0 \\ 0 & 2I
& I\\ 0 & (\hbar^2-1)I & 0\end{pmatrix}.
\end{eqnarray*}
If we identify $H$ with the closure span of $e'_1, \eta'_1, e'_2,
\eta'_2, ...$, and hence identify $H$ with $H\oplus H$, we see
that $Z(\hbar)^{-1}R(\hbar)Z(\hbar)$ is of the form
$\begin{pmatrix} A & W \\ B & 0\end{pmatrix}$,  where $A, W$ don't
depend on $\hbar$ and $W$ is an isometry with
$dim(ker(W^*))=\infty$  and $B$ is continuous for $0<\hbar\le 1$
with $\pa B\pa \le 1$.  By Lemma 5.1 of \cite{BP65} we can find
$X\in B(H)$ such that $A+WX=[B_1, B_2]$ for some $B_1, B_2\in
B(H)$ and $XW=0$. Replacing $B_1$ by some $B_1+\lambda I$ we may
assume that $B_1-I$ is invertible.  By the similarity
transformation
\begin{eqnarray*}
\begin{pmatrix} I & 0 \\ -X & I\end{pmatrix} \begin{pmatrix} A & W \\ B & 0\end{pmatrix} \begin{pmatrix} I & 0 \\ X & I\end{pmatrix}
=\begin{pmatrix} A+WX & W \\ -XA-XWX+B & 0\end{pmatrix}
\end{eqnarray*}
it suffices to find continuous paths $T(\hbar), S(\hbar)$ such that
\begin{eqnarray*}
[T(\hbar), S(\hbar)]=
\begin{pmatrix} A+WX & W \\ -XA-XWX+B & 0\end{pmatrix}
\end{eqnarray*}
and $\pa T(\hbar)\pa,\, \pa S(\hbar)\pa$ are bounded uniformly in $\hbar$.
Simple calculation shows that
\begin{eqnarray*}
T(\hbar)=\begin{pmatrix} B_1 & 0 \\ 0 & I\end{pmatrix},
S(\hbar)=\begin{pmatrix} B_2 & (B_1-I)^{-1}W \\
(-XA-XWX+B)(I-B_1)^{-1} & 0\end{pmatrix}
\end{eqnarray*}
satisfy the requirements.
\end{proof}

\begin{lemma} \label{asym rep 2:lemma}
Let $H$ be a separable Hilbert space. Then for any $n\in \Ne$
there exist norm continuous paths  $T_1(\hbar), \, T_2(\hbar),\,
\cdots, \, T_n(\hbar),\, S_1(\hbar),\, S_2(\hbar),\, \cdots,\, S_n(\hbar)$ of
operator in $B(H)$  for $0<\hbar\le 1$ such that

(1) $\lim_{\hbar\to 0}\pa [T_j(\hbar), S_k(\hbar)]-\delta^j_ki\pa
=0$,

(2) $[T_j(\hbar), T_k(\hbar)]=[S_j(\hbar), S_k(\hbar)]=0$,

(3) $\hbar^{-1/3}\pa T_j(\hbar)\pa$ and $\hbar^{-1/3}\pa
S_j(\hbar)\pa$ are bounded uniformly in $\hbar$.
\end{lemma}
\begin{proof}
Let $T(\hbar)$ and $S(\hbar)$ be as in Lemma~\ref{asym rep
1:lemma}.
Let
\begin{eqnarray*}
T_j(\hbar)&=&I\otimes I\otimes \cdots \otimes I\otimes
T(\hbar^{1/3})\otimes I\otimes \cdots \otimes I,\\
S_j(\hbar)&=&I\otimes I\otimes \cdots \otimes I\otimes
S(\hbar^{1/3})\otimes I\otimes \cdots \otimes I,
\end{eqnarray*}
 where
$T(\hbar^{1/3}), S(\hbar^{1/3})$ are at the j-th place. Identify $H$ with $H^{\otimes n}$.
Then clearly $T_j(\hbar), S_j(\hbar)$ satisfy the conditions.
\end{proof}

\begin{proof}[Proof of Theorem~\ref{asymptotic rep exit:thm}]
Let $H$ and $T_j(\hbar),\, S_j(\hbar)$ be as in Lemma~\ref{asym
rep 2:lemma}. For each $0<\hbar\le 1$ define a $\Re$-linear map
$\phi_{\hbar}:\Re^{2n}\rightarrow B(H)$ by
$\phi_{\hbar}(e_j)=T_j(\hbar), \phi_{\hbar}(e_{j+n})=S_j(\hbar)$
for $1\le j\le n$, where $e_1, \cdots, e_{2n}$ is the standard
basis of $\Re^{2n}$. Clearly $\{\phi_{\hbar}\}$ satisfy the
conditions (1)-(3) in Definition~\ref{Heisenberg:def}.

Denote the action of $U(n)$ on $\Re^{2n}$ by $\sigma$. Consider
the product $C^*$-algebra $\prod_{h\in U(n)}B(H)$, whose elements are
bounded maps $f:U(n)\rightarrow B(H)$. There is a natural
(discontinuous) action $\rho$ of $U(n)$ on $\prod_{h\in U(n)}B(H)$
given by $\rho_g(f)(h)=f(g^{-1}h)$. For each $0<\hbar\le 1$ define
a $\Re$-linear map $\varphi_{\hbar}:\Re^{2n}\rightarrow
\prod_{h\in U(n)}B(H)$ by
$\varphi_{\hbar}(u)(h)=\phi_{\hbar}(\sigma_{h^{-1}}(u))$ for $u\in
\Re^{2n}$. Clearly $\varphi_{\hbar}$ is $U(n)$-equivariant. Let
$\mathfrak{A}_{2n}$ be the $C^*$-subalgebra of $\prod_{h\in
U(n)}B(H)$ generated by $\cup_{0<\hbar\le
1}\varphi_{\hbar}(\Re^{2n})$. Then the restriction of $\rho$ on
$\mathfrak{A}_{2n}$ is continuous. Clearly $\pa
\varphi_{\hbar}\pa=\pa \phi_{\hbar}\pa$ and $\pa
\varphi_{\hbar}-\varphi_{\hbar'}\pa=\pa
\phi_{\hbar}-\phi_{\hbar'}\pa$, which verify the conditions (2)
and (3) of Definition~\ref{Heisenberg:def}. Using Lemma~\ref{asym
rep 2:lemma}(2) simple calculation shows that
\begin{eqnarray*}
\max_{1\le j<k\le 2n}\pa
[\varphi_{\hbar}(e_j),\varphi_{\hbar}(e_k)]-\delta^j_{k-n}i\pa \le
\max_{1\le j, k\le n}\pa [T_j(\hbar), S_k(\hbar)]-\delta^j_ki\pa.
\end{eqnarray*}
Then Definition~\ref{Heisenberg:def}(1) is also satisfied. In
particular, when $\hbar$ is small enough we have $\pa
[\varphi_{\hbar}(e_1), \varphi_{\hbar}(e_{n+1})]-i\pa<1$. Then
$[\varphi_{\hbar}(e_1), \varphi_{\hbar}(e_{n+1})]$ is invertible
in $\prod_{h\in U(n)}B(H)$ and hence $\mathfrak{A}_{2n}$ contain
the identity of $\prod_{h\in U(n)}B(H)$.
\end{proof}


\end{document}